\documentclass[11pt]{amsart}
\usepackage{latexsym}
\usepackage{amsmath}
\usepackage{amssymb}

\newcommand{\eproof}{\mbox{\ }\hfill $\Box$ \par \vskip 10pt}

\newtheorem{Theorem}{Theorem}[section]
\newtheorem{lemma}[Theorem]{Lemma}

\def\C{{\mathbb C}}
\def\R{{\mathbb R}}
\def\Re{{\rm Re}\:}
\def\Im{{\rm Im}\:}
\def\cal{\mathcal}

\begin{document}

\title{Localization of the interior transmission eigenvalues for a ball}
\author[V. Petkov]{Vesselin Petkov}
\author[G. Vodev]{Georgi Vodev}

\address {Universit\'e Bordeaux, Institut de Math\'ematiques de Bordeaux,  351, Cours de la Lib\'eration, 33405  Talence, France}
\email{petkov@math.u-bordeaux.fr}
\address {Universit\'e de Nantes, Laboratoire de Math\'ematiques Jean Leray, 2 rue de la Houssini\`ere, BP 92208, 44322 Nantes Cedex 03, France}
\email{Georgi.Vodev@univ-nantes.fr}

\date{}

\maketitle

\begin{abstract} 

We study the localization of the interior transmission eigenvalues (ITEs) in the case when the domain
is the unit ball $\{x \in \R^d:\: |x| \leq 1\}, \: d\geq 2,$ and the coefficients $c_j(x), \: j =1,2,$ and the indices of refraction $n_j(x), \: j =1,2,$ are constants near the boundary $|x| = 1$. We prove that in this case the eigenvalue-free
region obtained in \cite{kn:V2} for strictly concave domains can be significantly improved. In particular, if $c_j(x), n_j(x), j = 1,2$ are constants for $|x| \leq 1$,  we show that all (ITEs) lie in a strip $|\Im \lambda| \leq C$. 
\end{abstract}
\vspace{0.4cm}


\setcounter{section}{0}
\section{Introduction and statement of the result}

Let $\Omega\subset \R^d$, $d\ge 2$, be a bounded, connected domain with a $C^\infty$ smooth boundary $\Gamma=\partial\Omega$. 
A complex number $\lambda\neq 0$ with ${\rm Re}\,\lambda\ge 0$ will be called interior transmission eigenvalue (ITE) if the following problem has a non-trivial solution:
$$\left\{
\begin{array}{lll}
\left(\nabla c_1(x)\nabla+\lambda^2 n_1(x)\right)u_1=0 \:\mbox{in}\:\Omega,\\
\left(\nabla c_2(x)\nabla+\lambda^2 n_2(x)\right)u_2=0 \:\mbox{in}\:\Omega,\\
u_1=u_2,\,\,\, c_1\partial_\nu u_1=c_2\partial_\nu u_2\:\mbox{on}\:\Gamma,
\end{array}
\right.
\eqno{(1.1)}
$$
where $\nu$ denotes the Euclidean unit inner normal to $\Gamma$ and  $c_j(x) ,n_j(x) \in C^\infty(\overline\Omega)$, $j=1,2,$ are strictly positive real-valued functions.

The (ITEs) were first studied by Kirsch \cite{kn:K} and by Colton and Monk \cite{kn:CM} in the context of   
inverse scattering problems. It was  shown that the real (ITEs) correspond
to the frequencies for which the reconstruction algorithm in inverse scattering based
on the so-called linear sampling methods breaks down. This subject attracted the attention of many researchers and the number of publications devoted to the (ITEs) considerably increased in the recent ten years. The reader may consult the survey \cite{kn:CH} for a complete list of references and historical remarks.

It is well-known (e.g. see \cite{kn:PV}) that there exists a closed non-symmetric operator, $A$, associated in a natural
way to the problem (1.1),
such that the possible (ITEs) can be considered as the eigenvalues of $A$. The analysis of the (ITEs) leads to the following three problems:\\
$\:$\\
(A) Prove the discreteness of the spectrum of $A$ in $\C$;\\
(B) Find eigenvalue-free regions in $\C$;\\
(C)  Establish a Weyl formula for the counting function of all (ITEs)
$$N(r) = \#\{\lambda_j \:{\rm is}\: {\rm (ITE)},\,\:|\lambda_j | \leq r\}.$$
Note that the problem (A) is now relatively well studied (see \cite{kn:LC}, \cite{kn:S}, \cite{kn:PS}, \cite{kn:CLM} and the references therein). 
In fact, the problem (A) is reduced to that one of showing that the resolvent of $A$ is meromorphic with residues of finite rank. On the other
hand, this is
true (see \cite{kn:PV}) if the inverse of the operator $T(\lambda)$ introduced in Section 4 is meromorphic. The latter fact can be proved 
if the parametrix of
the operator $T(\lambda)$ constructed in the deep elliptic zone is invertible.

The problems (B) and (C) are more difficult, and they  
are of some interest for the numerical analysis of the (ITEs). In this direction it is interesting to find an optimal eigenvalue-free 
region and a Weyl formula with optimal remainder (see \cite{kn:R1}, \cite{kn:F}, \cite{kn:R2}, \cite{kn:LV} and the references therein). 
In a recent work \cite{kn:PV} the authors showed that (B) and (C) are closely related each other, and a {\it larger} eigenvalue-free region leads to a Weyl asymptotics with a {\it smaller} remainder term. More precisely, we proved that the remainder in the Weyl formula is
${\mathcal O}_\varepsilon(r^{d - \kappa+\varepsilon})$, $\forall\,0<\varepsilon\ll 1$, where $0<\kappa\le 1$ is such that there are no
(ITEs) in 
$$\left\{\lambda\in \C:\,{\rm Re}\,\lambda>1,\,|{\rm Im}\,\lambda|\ge C\left({\rm Re}\,\lambda\right)^{1-\kappa
}\right\}.$$
We conjecture that the optimal value of $\kappa$ must be $\kappa = 1.$

The present paper is devoted to the problem (B). More precisely, we are interested in finding as small as possible neighborhoods of the real axis containing all (ITEs). The first result of this type was obtained in \cite{kn:H} assuming $n_1(x) > 1$ in $\bar{\Omega}$ and $n_2(x) \equiv 1$, $c_1(x) \equiv c_2(x) \equiv 1$.
For domains $\Omega$ with {\it arbitrary geometry}, it has been shown in \cite{kn:V1} that under the condition (isotropic case)
$$c_1(x)\equiv c_2(x),\quad \partial_\nu c_1(x)\equiv \partial_\nu c_2(x),\quad n_1(x)\neq n_2(x)\quad\mbox{on}\quad\Gamma,\eqno{(1.2)}$$
or the condition (anisotropic case)
$$(c_1(x)-c_2(x))(c_1(x)n_1(x)-c_2(x)n_2(x))<0\quad\mbox{on}\quad\Gamma,\eqno{(1.3)}$$
there are no (ITEs) in the region 
$$\left\{\lambda\in \C:\,{\rm Re}\,\lambda>1,\,|{\rm Im}\,\lambda|\ge C_\varepsilon\left({\rm Re}\,\lambda\right)^{\frac{1}{2}+\varepsilon
}\right\},\quad\forall\, 0<\varepsilon\ll 1.$$
The localization of the (ITEs) has been recently studied in \cite{kn:V2} in the case when
the boundary $\Gamma$ is strictly concave with respect to both Riemannian metrics $\sum_{k=1}^d\frac{n_j(x)}{c_j(x)}dx_k^2$, $j=1,2$.  Under the conditions (1.2) or (1.3) it has been proved in \cite{kn:V2} that there are no (ITEs) in the region 
$$
\left\{\lambda\in \C:\,{\rm Re}\,\lambda>1,\,|{\rm Im}\,\lambda|\ge C_\varepsilon\left({\rm Re}\,\lambda\right)^\varepsilon
\right\},\quad\forall\, 0<\varepsilon\ll 1. \eqno{(1.4)}
$$
The approach in \cite{kn:V1} and \cite{kn:V2} is based on the construction of a semi-classical parametrix near the boundary for the problem
$$ \begin{cases} (h^2\nabla c(x)\nabla +zn(x)) u =0 \quad {\rm in }\quad \Omega,\\
u = f \quad {\rm on}\quad \Gamma, \end{cases} \eqno{(1.5)}$$
where $0 < h \ll 1$ is a semi-classical parameter and $z\in\C$ with $\Re z = 1$.
For domains with arbitrary geometry the parametrix construction for (1.5) works for $|\Im z| \geq h^{1/2-\epsilon}$, $0<\epsilon\ll 1$, while for
strictly concave domains, by a more complicated construction, one can cover the region $|\Im z| \geq h^{1-\epsilon}$.
It is a challenging problem to construct a semi-classical parametrix for (1.5) when $|\Im z | \geq Ch$, $C\gg 1$ being a constant.

The purpose of the present paper is to improve the eigenvalue-free region (1.4) in the case when the domain is the unit ball in $\R^d, d \geq 2.$ Given a parameter
$0<\delta\ll 1$, denote $\Omega(\delta)=\{x\in\overline{\Omega}:{\rm dist}(x,\Gamma)\le\delta\}$. Our main result is the following

\begin{Theorem} Let $\Omega=\{x\in \R^d:\,|x|\le 1\},\: d \geq 2,$ and suppose that there is a constant $0<\delta_0\ll 1$ such that the functions
$c_j(x)$, $n_j(x) $, $j=1,2$, are constants in $\Omega(\delta_0)$. Assume also either the condition $(1.2)$
or the condition $(1.3).$
Then, there is a constant $C>0$ such that there are no $(ITEs)$ in the region
$$
\{\lambda \in \C:\,{\rm Re}\,\lambda>1,\, |\Im \lambda| \geq C\log\left(\Re \lambda+1\right)\}.
\eqno{(1.6)}$$
If in addition the functions $c_j$, $n_j$, $j=1,2$, are constants everywhere in $\overline{\Omega}$, then 
there are no $(ITEs)$ in a larger region of the form
$$
\left\{\lambda\in \C:\,{\rm Re}\,\lambda>1,\,|{\rm Im}\,\lambda|\ge C\right\}.\eqno{(1.7)}
$$
\end{Theorem}
 
\noindent
{\bf Remark 1.} The eigenvalue-free region (1.6) is still valid if we add a compact cavity $K\subset\Omega$ and consider the equation
(1.1) in $\Omega\setminus K$ with Dirichlet condition on $\partial K$. Indeed, the only fact we need is the coercivity of the
corresponding Dirichlet realization (see the operator $G_D$ in Section 3), and this is used only in the proof of Lemma 3.4 below.\\

\noindent
{\bf Remark 2.} It is clear from the proof that the fact that the boundary $\Gamma$ is a sphere is not essential. In other words,
the eigenvalue-free regions (1.6) and (1.7) are still valid for any Riemannian manifold $\Omega=(0,1)\times\Gamma$ with metric
$g=dr^2+r^2\sigma$, where $r\in (0,1)$, and $(\Gamma,\sigma)$ is an arbitrary $(d-1)$-dimensional Riemannian manifold without boundary, the metric $\sigma$ being independent of $r$.\\

\noindent
In the isotropic case when $c_j \equiv 1, j = 1, 2$  and $n_1 = 1, n_2 \neq 1$ is constant, the eigenvalue-free region (1.7) 
has been established in the one-dimensional case $\Omega = \{x \in \R: |x| \leq 1\}$ (see \cite{kn:S}, \cite{kn:PS}). 
Moreover, the case of the ball $\{x \in \R^d:\: |x| \leq 1$\}, $d = 2, 3,$ and radial refraction indices have been studied in 
\cite{kn:LC}, \cite{kn:CL}, \cite{kn:CLM}, where spherical symmetric eigenfunctions depending only on 
the radial variable $r = |x|$ has been considered. For example, the analysis of such eigenfunctions in $\R^3$ leads to 
the following one-dimensional problem 
$$\begin{cases} u'' + \frac{2}{r} u'+ \lambda^2 n(r) u = 0 ,\quad 0 < r < 1,\\
v'' + \frac{2}{r} v'+ \lambda^2 v = 0,\quad 0 < r < 1,\\
u(1) = v(1),\: u'(1) = v'(1),\end{cases}\eqno{(1.8)}$$
where $n(r)$ is a strictly positive function. Among other things, it was shown in \cite{kn:CLM} that if $n(1)=1$ and $n'(1)$ or 
$n''(1)$ is non-zero, then there may exist
infinitely many complex eigenvalues of the problem (1.8) lying outside any strip parallel to the real axis. This example shows that in
the isotropic case 
the condition $n(1)\neq 1$ (resp. (1.2)) is important to have an eigenvalue-free region like (1.7).
It also follows from the analysis in \cite{kn:LC} (see Sections 3 and 4) that when $n=Const$ and $\sqrt{n}$ is a rational number belonging
to the interval $(1,2)$, then there
exists a sequence of (ITEs), $\lambda_k=\alpha k+\beta$, $k=0,1,2,...,$ with some constants $\alpha>0$ and $\beta\in\C$,
${\rm Im}\,\beta\neq 0$. This example shows that
the eigenvalue-free region (1.7) is sharp and cannot be improved in general.

To study all (ITEs) and all eigenfunctions, however, one has to consider a family of infinitely many one-dimensional problems. 
Such an analysis is carried out in
\cite{kn:PS} in the isotropic case when the domain is the ball $\{x \in \R^d:\: |x| \leq 1\},\, d\geq 1,$ and 
$$c_1\equiv c_2 \equiv 1, n_1 \equiv 1, n_2 \equiv \gamma^2=Const \neq 1.$$ 
In this case the (ITEs) are the zeros in $\C$ of the family of functions
$$F_{\nu}(\lambda) = \gamma J_{\nu}(\lambda)J'_{\nu}(\gamma \lambda) - J_{\nu}(\gamma\lambda) J'_{\nu}(\lambda),
\quad \nu = l + d/2 - 1,\quad l=0,1,2,...,$$
where $J_\nu$ denotes the Bessel function of order $\nu$. It has been proved in \cite{kn:PS} that there are infinitely many real 
(ITEs) whose counting function has a Weyl asymptotics. When $d=1$ a Weyl asymptotics for the counting function of all
(ITEs) is also obtained.

 To prove  Theorem 1.1 we follow the same strategy as in \cite{kn:V1}, \cite{kn:V2}, which consists of deriving the eigenvalue-free
 region from some approximation properties of the {\it interior} Dirichlet-to-Neumann (DN) map. In our case we have to
 approximate the DN map 
 $${\cal N}_0(\lambda):H^{s+1}(\Gamma)\to H^s(\Gamma)$$
  for the domain $\Omega=\{x\in \R^d:\,|x|\le 1\}$ defined by
  $${\cal N}_0(\lambda)f:=\lambda^{-1}\partial_\nu u|_{\Gamma}$$
  where $\nu$ is the unit inner normal to $\Gamma=\partial\Omega$ and $u$ solves the equation
$$\left\{
\begin{array}{lll}
 \left(\Delta+\lambda^2\right)u=0&\mbox{in}& \Omega,\\
 u=f&\mbox{on}&\Gamma,
\end{array}
\right.
$$
$\Delta$ being the negative Euclidean Laplacian. 
 Recall that the interior DN map is 
 a meromorphic operator-valued function with poles lying on the positive real axis.
 Thus, the eigenvalue-free region turns out to be the region in which the DN map can be approximated by a simpler operator
 of the form $f(\Delta_\Gamma)$, where $f$ is a complex-valued function and $\Delta_\Gamma$ denotes the negative Laplace-Beltrami operator on the boundary $\Gamma$ equipped with the Riemannian metric induced by the Euclidean one. With such an approximation in hands,
 the problem of proving the eigenvalue-free region is transformed into the much simpler one of inverting complex-valued functions, which in turn is done using the conditions (1.2) or (1.3) (see Section 4). Therefore, a large portion of the present paper is devoted to the study of 
 the interior DN map ${\cal N}_0(\lambda)$ using the Bessel functions. Thus, instead of a parametrix we have an exact formula of the DN map
 (see Theorem 3.1 and its proof). Then we use the asymptotic expansions of the Bessel functions in terms of the Airy function
 to get the desired approximation (see Theorems 2.1 and 3.1). Of course, we cannot proceed in this way when the coefficients
 are supposed to be constants only in a neighborhood of the boundary. In this latter case we show that the DN map can be approximated by
 the DN map associated to the corresponding problem with constant coefficients everywhere and for which we have an
 explicit expression in terms of the Bessel functions (see Lemma 3.4). 
 
 We expect that the eigenvalue-free regions (1.6) and (1.7) are still true for more general domains, but this is hard to prove
 because the available semi-classical parametrix constructions for the DN map lead to the existence of smaller regions (see 
 \cite{kn:V1}, \cite{kn:V2}).

 \section{Some properties of the Bessel functions}
 
 We begin this section by recalling some basic properties of the Bessel functions $J_\nu(z)$ of real order $\nu\ge 0$ (e.g. see \cite{kn:O}). 
 The function $J_\nu(z)$ satisfies the equation
 $$\left(\partial_z^2+z^{-1}\partial_z+1-(\nu/z)^2\right)u(z)=0.$$
 Then the function $b_\nu(z)=z^{1/2}J_\nu(z)$ satisfies the equation
 $$\partial_z^2v+\left(1-\frac{\nu^2-1/4}{z^2}\right)v=0.$$
  For $z\in \C$ with $\Re z>0$, $\Im z\neq 0$, and an integer $k\ge 0$, set
  $$\psi_\nu(z)=\frac{J'_\nu(z)}{J_\nu(z)},\quad \eta^{(k)}_{\nu}(z)=\frac{J^{(k)}_\nu(\kappa z)}{J_\nu(z)},$$
 where $J^{(k)}_\nu(z)=\frac{d^kJ_\nu(z)}{dz^k}$ and $0<\kappa<1$ is a parameter independent of $z$ and $\nu$. Clearly, the functions 
 $\eta^{(k)}_{\nu}(z)$ 
depend on $\kappa$ but for simplicity of the notations we will omit to note this.  Set also $\rho(z)=\sqrt{z^2-1}$ with ${\rm Re}\,\rho>0$. Our goal in this section is to prove the following
 
 \begin{Theorem} For every $0<\delta\ll 1$,  
 there are positive constants $C_\delta$, $C'_\delta$ and $\delta_1$ 
 such that for ${\rm Re}\,\lambda\ge C_\delta$, $C'_\delta\le|{\rm Im}\,\lambda|\le\delta_1{\rm Re}\,\lambda$, 
 $\nu\ge 0$, we have the estimate
 $$(1+\nu/|\lambda|)\left|\psi_{\nu}(\lambda)-\rho(\nu/\lambda)\right|\le \delta.\eqno{(2.1)}$$
 There exist also positive constants $C$, $C'$, $C_1$, $C_2$ and $\delta_1$ independent of  $\nu$ but depending on $\kappa$
 such that for ${\rm Re}\,\lambda\ge C_1$, $C_2\le|{\rm Im}\,\lambda|\le\delta_1{\rm Re}\,\lambda$, 
 $\nu\ge 0$, we have the estimate
 $$(1+\nu/|\lambda|)^2|\eta^{(0)}_{\nu}(\lambda)|+|\eta^{(1)}_{\nu}(\lambda)|+|\eta^{(2)}_{\nu}(\lambda)|\le C'|\lambda|^{1/3} e^{-C|{\rm Im}\,\lambda|}.\eqno{(2.2)}$$
 \end{Theorem}
 
 {\it Proof.} We will consider several cases.
 
 {\bf Case 1.} $0\le\nu\le C_0$ with some constant $C_0 > 0$. We have $2J_\nu(\lambda)=H^+_\nu(\lambda)+H^-_\nu(\lambda)$, where $H^\pm_\nu(\lambda)$ are the Hankel functions \footnote{$H^{\pm}_{\nu}(\lambda)$ are the Hankel functions of first and second kind} having the asymptotic expansions
 (see (4.03) and (4.04), p.238 in \cite{kn:O})
 $$H^\pm_\nu(\lambda)=\left(\frac{2}{\pi\lambda}\right)^{1/2}e^{\pm i\lambda}q_\nu^\pm(\lambda),$$
 $$q_\nu^\pm(\lambda)=e^{\pm i(-\nu\pi/2-\pi/4)}\sum_{s=0}^\infty\left(\frac{\pm i}{\lambda}\right)^sA_s(\nu),\eqno{(2.3)}$$
 where all $A_s(\nu)$ are real, $A_0(\nu)=1$, $A_1(\nu) =\frac{4\nu^2-1}{8}$. Moreover, $q_{1/2}^\pm(\lambda)=\pm i$.
 All derivatives of $q_\nu^\pm(\nu)$ have asymptotic expansions obtained by differentiating (2.3). 
 Without loss of generality, we may suppose that ${\rm Im}\,\lambda>0$. For $\nu\neq 1/2$ we have
 $$\left|\frac{q_\nu^+(\lambda)}{q_\nu^-(\lambda)}\right|= 1+{\cal O}(|\lambda|^{-1}),\quad
 \left|\frac{(q_\nu^+)'(\lambda)}{(q_\nu^-)'(\lambda)}\right|= 1+{\cal O}(|\lambda|^{-1}),\quad
 \left|\frac{(q_\nu^-)'(\lambda)}{(q_\nu^-)(\lambda)}\right|= {\cal O}(|\lambda|^{-2}),$$
 $$\left|1 +e^{2i\lambda}\frac{q_\nu^+(\lambda)}{q_\nu^-(\lambda)}\right|\ge 1-e^{-2{\rm Im}\,\lambda}
 \left|\frac{q_\nu^+(\lambda)}{q_\nu^-(\lambda)}\right|\ge 1-\frac{3}{2}e^{-2{\rm Im}\,\lambda}\ge \frac{1}{2}$$
 provided $|\lambda|$ and ${\rm Im}\,\lambda$ are taken large enough.
 We can write the function $\psi_\nu$ as follows
 $$\psi_\nu(\lambda)+(2\lambda)^{-1}=i\frac{e^{i\lambda}q_\nu^+(\lambda)-e^{-i\lambda}q_\nu^-(\lambda)}{e^{i\lambda}q_\nu^+(\lambda)
 +e^{-i\lambda}q_\nu^-(\lambda)}+\frac{e^{i\lambda}(q_\nu^+)'(\lambda)+
 e^{-i\lambda}(q_\nu^-)'(\lambda)}{e^{i\lambda}q_\nu^+(\lambda)+e^{-i\lambda}q_\nu^-(\lambda)}.$$
 By using the above inequalities, we get
 $$|\psi_\nu(\lambda)+i|\le C|\lambda|^{-1}+Ce^{-2{\rm Im}\,\lambda}.\eqno{(2.4)}$$
 Since in this case $\rho(\nu/\lambda)=-i+{\cal O}(|\lambda|^{-2})$, the estimate (2.1) follows from (2.4). The estimate (2.2) for
 $|\eta^{(0)}_{\nu}(\lambda)|$ follows
 in the same way from the formula
 $$\frac{J_\nu(\kappa\lambda)}{J_\nu(\lambda)}=\kappa^{-1/2}\frac{e^{i\kappa\lambda}q_\nu^+(\kappa\lambda)
 +e^{-i\kappa\lambda}q_\nu^-(\kappa\lambda)}{e^{i\lambda}q_\nu^+(\lambda)
 +e^{-i\lambda}q_\nu^-(\lambda)}.$$
 Indeed, as above, one can easily see that $\eta_{\nu}^{(0)}(\lambda)={\cal O}\left(e^{-(1-\kappa){\rm Im}\,\lambda}\right)$.
 
 {\bf Case 2.} $\nu\gg 1$. We set $z=\lambda/\nu$. Then $1/\nu\ll|{\rm Im}\,z|\ll{\rm Re}\,z$. 
 In this case we will use the asymptotic expansions of the Bessel functions in terms of the Airy function
 ${\rm Ai}(\sigma)$. Recall first that ${\rm Ai}(\sigma)$ has the expansion
 $${\rm Ai}(\sigma)=\sigma^{-1/4}e^{-\frac{2}{3}\sigma^{3/2}}\sum_{\ell=0}^\infty\beta_\ell\,\sigma^{-3\ell/2}\eqno{(2.5)}$$
 for $|\sigma|\gg 1$, $\sigma\in\Lambda_\varepsilon:=\{\sigma\in \C:|{\rm arg}\,\sigma|\le \pi-\varepsilon\}$,
 $0<\varepsilon\ll 1$, where $\beta_\ell$ are real numbers and the fractional powers of $\sigma$ take their principal values. 
 The expansion (2.5) implies
 $$F(\sigma):=\frac{{\rm Ai}'(\sigma)}{{\rm Ai}(\sigma)}=-\sigma^{1/2}\sum_{\ell=0}^\infty\widetilde\beta_\ell\,\sigma^{-3\ell/2}\eqno{(2.6)}$$
 where $\widetilde\beta_0=1$, $\widetilde\beta_1=1/4$. The behavior of the function $F$ in $\C \setminus\Lambda_\varepsilon$ is more complicated and is given by the following
 
  \begin{lemma}
 For $\sigma\in \C\setminus\Lambda_\varepsilon$, $\Im \sigma\neq 0$,
  we have the bounds 
 $$|F(\sigma)|\le C|\sigma|^{1/2}+C|{\rm Im}\,\sigma|^{-1},\eqno{(2.7)}$$
 $$|{\rm Ai}(\sigma)|\le C\langle\sigma\rangle^{-1/4}e^{\frac{2}{3}|{\rm Re}\,\sigma^{3/2}|},\eqno{(2.8)}$$
 $$|{\rm Ai}(\sigma)|^{-1}\le C\langle\sigma\rangle^{-1/4}\left(|\sigma|^{1/2}+
 |{\rm Im}\,\sigma|^{-1}\right)e^{-\frac{2}{3}|{\rm Re}\,\sigma^{3/2}|},\eqno{(2.9)}$$
 where we have used the notation $\langle\sigma\rangle=(1+|\sigma|^2)^{1/2}$. 
 For $\sigma\in \C \setminus\Lambda_\varepsilon$, $|\sigma|\gg 1$, $|{\rm Re}\,\sigma^{3/2}|\gg 1$,
  we have the bound 
 $$\left|F(\sigma)+\sigma^{1/2}+\frac{1}{4\sigma}\right|\le C|\sigma|^{1/2}e^{-|{\rm Re}\,\sigma^{3/2}|}.\eqno{(2.10)}$$
 \end{lemma}
 
 {\it Proof.} The bound (2.7) is proved in \cite{kn:V2} (see Lemma 3.1). 
 To prove the other bounds, we will use that ${\rm Ai}(-\sigma)=
 {\rm Ai}_+(\sigma)+{\rm Ai}_-(\sigma)$, where ${\rm Ai}_\pm(\sigma)=e^{\pm \pi i/3}{\rm Ai}\left(\sigma e^{\pm\pi i/3}\right)$.
 By (2.5), for $|\arg\sigma|\le\varepsilon$, $|\sigma|\gg 1$, we have
 $${\rm Ai}_\pm(\sigma)=\sigma^{-1/4}e^{\pm i\frac{2}{3}\sigma^{3/2}}a_\pm(\sigma),\quad a_\pm(\sigma)=\sum_{\ell=0}^\infty \beta_\ell^\pm\,\sigma^{-3\ell/2}\eqno{(2.11)}$$
 with $|\beta_\ell^\pm|=|\beta_\ell|$. In particular, this implies 
 $$|{\rm Ai}_\pm(\sigma)|\le C\langle\sigma\rangle^{-1/4}e^{\mp\frac{2}{3}{\rm Im}\,\sigma^{3/2}},\quad |{\rm Ai}'_\pm(\sigma)|\le C
 \langle\sigma\rangle^{1/4}e^{\mp\frac{2}{3}{\rm Im}\,\sigma^{3/2}}.\eqno{(2.12)}$$
 Since $|{\rm Im}\,\sigma^{3/2}|=|{\rm Re}\,(-\sigma)^{3/2}|$, 
 we get (2.8) from (2.12). The bound (2.9) follows from (2.7), (2.12) and the identity
 $${\rm Ai}(-\sigma)^{-1}=c_\pm F(-\sigma){\rm Ai}_\pm(\sigma)+\widetilde c_\pm{\rm Ai}'_\pm(\sigma) \eqno{(2.13)}$$
 where $c_\pm$ and $\widetilde c_\pm$ are some constants. To prove (2.10), observe that, if 
 $|\arg\sigma|\le\varepsilon$, ${\rm Im}\,\sigma>0$, 
 we have ${\rm Im}\,\sigma^{3/2}>0$, and we can write
 $$-F(-\sigma)+i\sigma^{1/2}+\frac{1}{4\sigma}=2i\sigma^{1/2}\frac{e^{ i\frac{2}{3}\sigma^{3/2}}a_+(\sigma)}{e^{ i\frac{2}{3}\sigma^{3/2}}a_+(\sigma)+e^{-i\frac{2}{3}\sigma^{3/2}}a_-(\sigma)}$$
 $$+\frac{e^{ i\frac{2}{3}\sigma^{3/2}}a'_+(\sigma)+e^{- i\frac{2}{3}\sigma^{3/2}}a'_-(\sigma)}{e^{ i\frac{2}{3}\sigma^{3/2}}a_+(\sigma)+e^{-i\frac{2}{3}\sigma^{3/2}}a_-(\sigma)}.\eqno{(2.14)}$$
 The above expansions imply
 $$\left|\frac{a_-(\sigma)}{a_+(\sigma)}\right|= 1+{\cal O}(|\sigma|^{-1}),\quad
 \left|\frac{a'_-(\sigma)}{a'_+(\sigma)}\right|= 1+{\cal O}(|\sigma|^{-1}),\quad
 \left|\frac{a'_+(\sigma)}{a_+(\sigma)}\right|= {\cal O}(|\sigma|^{-1}).$$
 Therefore in this case (2.10) follows from (2.14) after making a change of variables $\sigma\to -\sigma$
 and using that if $|\arg\sigma|\le\varepsilon$, ${\rm Im}\,\sigma>0$, then $-\sigma\in\C\setminus\Lambda_\varepsilon$ and 
 $(-\sigma)^{1/2}=-i\sigma^{1/2}$. The analysis of the case ${\rm Im}\,\sigma<0$ is similar.
 \eproof
 
 Define the functions $\varphi(z)$ and $\zeta(z)$ by
 $$\varphi=\frac{2}{3}\zeta^{3/2}=\ln \frac{1+(1-z^2)^{1/2}}{z}-(1-z^2)^{1/2},\quad\quad |\arg z|<\pi,$$
 where the branches take their principal values when $z\in(0,1)$, $\varphi,\zeta\in (0,+\infty)$, and are continuous elsewhere.
 It is well-known (e.g. see pages 420-422 of \cite{kn:O}) that the function $\zeta(z)$ is holomorphic for $|\arg z|<\pi$, $\zeta(z)$ takes real values for $z\in(0,+\infty)$, and $\zeta(z)=2^{1/3}(1-z)+{\cal O}(|1-z|^2)$ in a neighborhood of $z=1$.
 Moreover, $\zeta(z)\to -\infty$ as $z\to +\infty$ and $\zeta(z)\to +\infty$ as $z\to 0^+$. 
  The first derivatives of $\varphi(z)$ and $\zeta(z)$ satisfy
 $$\zeta(z)^{1/2}\zeta'(z)=\varphi'(z)=-\frac{(1-z^2)^{1/2}}{z}.\eqno{(2.15)}$$
 One can easily see that for $0<\pm{\rm Im}\,z\ll {\rm Re}\,z$
 we have
 $${\rm Re}\,\varphi'(z)<0,\quad \pm {\rm Im}\,\varphi'(z)>0.\eqno{(2.16)}$$
 In particular, this implies that the function $\rho$ defined above satisfies
 $$\rho\left(\frac{1}{z}\right)=\frac{(1-z^2)^{1/2}}{z}.\eqno{(2.17)}$$
 Given parameters $0<\delta,\delta_1\ll 1$, set $$\Theta_1(\delta,\delta_1)=\left\{{\rm Re}\,z\ge 1+\delta^2,\,0<|{\rm Im}\,z|\le\delta_1{\rm Re}\,z\right\},$$
 $$\Theta_2(\delta,\delta_1)=\left\{0<{\rm Re}\,z\le 1-\delta^2,\,0<|{\rm Im}\,z|\le\delta_1{\rm Re}\,z\right\},$$
  $$\Theta_0(\delta,\delta_1)=\left\{1-\delta^2\le{\rm Re}\,z\le 1+\delta^2,\,0<|{\rm Im}\,z|\le\delta_1{\rm Re}\,z\right\}.$$
  The next lemma is more or less well-known and follows from the properties of the functions
  $\varphi$ and $\zeta$ studied in \cite{kn:O}. We will sketch the proof for the sake of completeness.
  
  \begin{lemma} For every $0<\delta\ll 1$ there is $\delta_1=\delta_1(\delta)>0$ such that the following properties hold:
    For $z\in \Theta_1(\delta,\delta_1)$ we have $|\arg \zeta(z)|=\pi-{\cal O}(\delta)$, and
 $$2|\zeta(z)|^{1/2}|{\rm Im}\,\zeta(z)|\ge |{\rm Re}\,\varphi(z)|\ge C|{\rm Im}\,z|\eqno{(2.18)}$$
 with a constant $C>0$ depending on $\delta$.
For $z\in \Theta_2(\delta,\delta_1)$ we have $|\arg \zeta(z)|={\cal O}(\delta)$. 
For $z\in \Theta_0(\delta,\delta_1)$ we have
$$|{\rm Im}\,\zeta(z)|\ge  |{\rm Im}\,z|.\eqno{(2.19)}$$
 \end{lemma}
 
 {\it Proof.} We will use the formula
 $$\varphi(z)-\varphi({\rm Re}\,z)=\int_0^1\frac{d}{d\tau}\varphi({\rm Re}\,z+i\tau{\rm Im}\,z)d\tau
 =i{\rm Im}\,z\int_0^1\varphi'({\rm Re}\,z+i\tau{\rm Im}\,z)d\tau.\eqno{(2.20)}$$
 Let $z\in \Theta_1(\delta,\delta_1)$. Then 
 $${\rm Re}\,\varphi({\rm Re}\,z)=0,\quad {\rm Im}\,\varphi({\rm Re}\,z)\ge C_\delta{\rm Re}\,z.$$ 
 In this case we also have
 $$\varphi'({\rm Re}\,z+i\tau{\rm Im}\,z)={\cal O}_\delta(1)$$
 and, in view of (2.16), if $\pm{\rm Im}\,z>0$, 
 $$\pm {\rm Im}\,\varphi'({\rm Re}\,z+i\tau{\rm Im}\,z)\ge C_\delta-{\cal O}_\delta(\delta_1)\ge C_\delta/2>0$$
 provided $\delta_1$ is taken small enough. Thus, by (2.20) we get
 $$-{\rm Re}\,\varphi(z)\ge C_\delta|{\rm Im}\,z|,$$
 $$\pm{\rm Im}\,\varphi(z)\ge (C_\delta-{\cal O}_\delta(\delta_1)){\rm Re}\,z\ge 2^{-1}C_\delta{\rm Re}\,z,\quad \pm{\rm Im}\,z>0.$$ 
 This yields ${\rm Re}\,(\mp i\varphi(z))>0$, $\pm{\rm Im}\,(\mp i\varphi(z))>0$, and hence $0<\pm\arg (\mp i\varphi(z))={\cal O}_\delta(\delta_1)={\cal O}(\delta)$ if $\delta_1$ is small enough. Since
 $$\varphi=\frac{2}{3}\zeta^{3/2}=\pm i\frac{2}{3}(-\zeta)^{3/2},$$
 we have
 $$0<\pm\arg(-\zeta(z))=\frac{2}{3}\arg (\mp i\varphi(z))={\cal O}(\delta)$$
 and
 $$|{\rm Re}\,\varphi(z)|=\frac{2}{3}\left|{\rm Im}\,(-\zeta(z))^{3/2}\right|=|{\rm Im}\,\zeta(z)||\zeta(z)|^{1/2}(1+{\cal O}(\delta)).$$
 Let $z\in \Theta_2(\delta,\delta_1)$. Then
  $${\rm Im}\,\varphi({\rm Re}\,z)=0,\quad {\rm Re}\,\varphi({\rm Re}\,z)\ge C_\delta>0,$$ 
  $${\rm Im}\,\varphi'({\rm Re}\,z+i\tau{\rm Im}\,z)={\cal O}_\delta(1),$$
 $$-{\rm Re}\,\varphi'({\rm Re}\,z+i\tau{\rm Im}\,z)\ge (C_\delta-{\cal O}_\delta(\delta_1))({\rm Re}\,z)^{-1}
 \ge 2^{-1}C_\delta({\rm Re}\,z)^{-1},$$
 provided $\delta_1$ is taken small enough. Therefore, by (2.20) we get
 $$|{\rm Im}\,\varphi(z)|\le C_\delta\frac{|{\rm Im}\,z|}{{\rm Re}\,z}={\cal O}_\delta(\delta_1),$$
 $${\rm Re}\,\varphi(z)={\rm Re}\,\varphi({\rm Re}\,z)+{\cal O}_\delta(|{\rm Im}\,z|)\ge C_\delta-{\cal O}_\delta(\delta_1)\ge C_\delta/2.$$
 Hence, $\arg\varphi(z)={\cal O}_\delta(\delta_1)={\cal O}(\delta)$, which yields
 $$\arg\zeta(z)=\frac{2}{3}\arg \varphi(z)={\cal O}(\delta).$$
 Let $z\in \Theta_0(\delta,\delta_1)$. Then we have $\zeta'(z)=-2^{1/3}+{\cal O}(|1-z|)$ at $z=1$.
 To prove (2.19) we will use the formula
 $$\zeta(z)-\zeta({\rm Re}\,z)=\int_0^1\frac{d}{d\tau}\zeta({\rm Re}\,z+i\tau{\rm Im}\,z)d\tau
 =i{\rm Im}\,z\int_0^1\zeta'({\rm Re}\,z+i\tau{\rm Im}\,z)d\tau$$
 $$=-i2^{1/3}{\rm Im}\,z(1+{\cal O}(\delta)).\eqno{(2.21)}$$
 Since ${\rm Im}\,\zeta({\rm Re}\,z)=0$, we deduce from (2.21),
 $${\rm Im}\,\zeta(z)=-2^{1/3}{\rm Im}\,z(1+{\cal O}(\delta))$$
 which clearly implies (2.19).
 \eproof
 
 For $|\arg z|\le\varepsilon$, $\nu\to +\infty$, we have the expansion (see \cite{kn:O}, (10.18), p.423 and more generally (9.02), p.418)
 $$J_\nu(\nu z)=2^{1/2}\nu^{-1/3}\left(\frac{\zeta}{1-z^2}\right)^{1/4}\left({\rm Ai}(\nu^{2/3}\zeta)
 A(\zeta)+\nu^{-4/3}{\rm Ai'}(\nu^{2/3}\zeta)B(\zeta)+{\cal E}_1(\nu, \zeta)\right)$$
 where
 $$A(\zeta)=\sum_{s=0}^M\frac{A_s(\zeta)}{\nu^{2s}}, \quad\quad B(\zeta)=\sum_{s=0}^M\frac{B_s(\zeta)}{\nu^{2s}},$$
 for every integer $M\gg 1$, 
where the functions $A_s(\zeta)$, $B_s(\zeta)$ are smooth and bounded with their derivatives, $A_0(\zeta)=1$,
 $B_s(\zeta)={\cal O}(\langle\zeta\rangle^{-1/2})$. The error term satisfies the bounds (see \cite{kn:O}, (10.19), p.423 and (9.03), p.418
 together with the notations on p.415)
 $$\left|\partial_\zeta^{\ell}{\cal E}_1(\nu, \zeta)\right|\le C_M\nu^{-2M}\langle\zeta\rangle^{(\ell-1)/4}
  e^{\frac{2\nu}{3}|{\rm Re}\,\varphi(z)|},\quad \ell=0,1.
 \eqno{(2.22)}$$
 We will derive now a similar expansion for the first derivative of $J_\nu$. To this end, observe first that by (2.15) we have
 $$\left(\frac{\zeta}{1-z^2}\right)^{1/4}\zeta'(z)=-\frac{1}{z}\left(\frac{\zeta}{1-z^2}\right)^{-1/4},$$
 $$\frac{\partial}{\partial z}\left(\frac{\zeta}{1-z^2}\right)^{1/4}=-\frac{1}{z}\left(\frac{\zeta}{1-z^2}\right)^{-1/4}\phi(z),$$
 where
 $$\phi(z)=\frac{1}{4\zeta}-\frac{\zeta^{1/2}z^2}{2(1-z^2)^{3/2}}.$$
 Since $|\zeta|\sim|z|$ as $|z|\to \infty$, $|\zeta|\sim\log(|z|^{-1})$ as $|z|\to 0$, $\zeta(z)=2^{1/3}(1-z)+{\cal O}(|1-z|^2)$ as $z\to 1$, we have
 $$\zeta^{-1/2}\left(\phi(z)-\frac{1}{4\zeta}\right)=\left\{
 \begin{array}{lll}
 {\cal O}_\epsilon(|z|^{2-\epsilon}),\,\forall 0<\epsilon\ll 1,\quad |z|\to 0,\\
 {\cal O}(\langle\zeta\rangle^{-1}),\quad |z|\to \infty,\\
 {\cal O}(|\zeta|^{-3/2}),\quad z\to 1.
\end{array}
\right.$$
 Differentiating the expansion of $J_\nu$ above with respect to the variable $z$ and using that ${\rm Ai}''(\sigma)=
 \sigma {\rm Ai}(\sigma)$, we get 
 $$z(J_\nu)'(\nu z)=-2^{1/2}\nu^{-2/3}\left(\frac{\zeta}{1-z^2}\right)^{-1/4}\left({\rm Ai}'(\nu^{2/3}\zeta)
 C(\zeta)+\nu^{-2/3}{\rm Ai}(\nu^{2/3}\zeta)D(\zeta)+{\cal E}_2(\nu, \zeta)\right)$$
 where
 $$C=A+\nu^{-2}(\partial_\zeta B+\phi B),\quad D=\partial_\zeta A+\phi A+\zeta B,\quad {\cal E}_2=\nu^{-2/3}(\partial_\zeta{\cal E}_1 +\phi
 {\cal E}_1).$$
 Then we have the identity
 $$\psi_\nu(\nu z)-\frac{(1-z^2)^{1/2}}{z}$$ $$=-\left(\frac{(1-z^2)^{1/2}}{z}\right)
 \frac{\Phi(\zeta)(1+P_1(\zeta))
 +P_2(\zeta)+P_3(\zeta)}{1+Q_1(\zeta)+\nu^{-1/3}\zeta^{-1/2}F(\nu^{2/3}\zeta)Q_2(\zeta)+Q_3(\zeta)}$$
 where
 $$\Phi(\zeta)=\nu^{-1/3}\zeta^{-1/2}F(\nu^{2/3}\zeta)+1+(4\nu\zeta^{3/2})^{-1},$$
 $$Q_1(\zeta)=A(\zeta)-1={\cal O}(\nu^{-2}),$$
 $$Q_2(\zeta)=\nu^{-1}\zeta^{1/2}B(\zeta)={\cal O}(\nu^{-1}w(\zeta)^{1/2}),$$
 $$Q_3(\zeta)={\cal E}_1(\nu, \zeta){\rm Ai}(\nu^{2/3}\zeta)^{-1},$$
 $$P_1(\zeta)=C(\zeta)-1+\nu^{-1}\zeta^{1/2}B(\zeta)$$
 $$=A(\zeta)-1+\nu^{-1}\zeta^{1/2}B(\zeta)+\nu^{-2}(\partial_\zeta B(\zeta)+\phi B(\zeta))={\cal O}(\nu^{-1}),$$
 $$P_2(\zeta)=\left(1+(4\nu\zeta^{3/2})^{-1}\right)(A-C)-(4\nu\zeta^{3/2})^{-1}A$$ $$-\left(1+(4\nu\zeta^{3/2})^{-1}\right)Q_2+
 \nu^{-1}\zeta^{-1/2}D$$
 $$=\nu^{-2}\left(1+(4\nu\zeta^{3/2})^{-1}\right)(\partial_\zeta B(\zeta)+\phi B(\zeta))-(4\nu\zeta^{3/2})^{-1}(A(\zeta)-1)$$
 $$-\nu^{-1}(4\nu\zeta^{3/2})^{-1}\zeta^{1/2}B(\zeta)+
 \nu^{-1}\zeta^{-1/2}(\partial_\zeta A(\zeta)+\phi (A(\zeta)-1))$$
 $$+\nu^{-1}\zeta^{-1/2}\left(\phi-(4\zeta)^{-1}\right)$$ $$={\cal O}\left(\nu^{-1}w(\zeta)^{-3/2}w(z)^{2-\epsilon}\right)+
 {\cal O}\left(\nu^{-2}\right),$$
 $$P_3(\zeta)=\nu^{-1/3}\left(\zeta^{-1/2}{\cal E}_2(\nu, \zeta)+{\cal E}_1(\nu, \zeta)\right){\rm Ai}(\nu^{2/3}\zeta)^{-1}$$
 $$=\nu^{-1}\left(\zeta^{-1/2}(\partial_\zeta{\cal E}_1(\nu, \zeta)+\phi
 {\cal E}_1(\nu, \zeta))+\nu^{1/3}{\cal E}_1(\nu, \zeta)\right){\rm Ai}(\nu^{2/3}\zeta)^{-1}$$
 uniformly for $|\zeta|\ge \nu^{-1}$, where $w(\sigma)=|\sigma|/\langle \sigma\rangle$. 
 We will consider now three cases.
 
 a) $z\in\Theta_1(\delta,\delta_1)$. Then $|\zeta|\ge C_\delta>0$, and by Lemma 2.3 we have $|\arg \zeta(z)|=\pi-{\cal O}(\delta)$ and  $|{\rm Im}\,z|\gg\nu^{-1}$ implies
  $\nu|{\rm Re}\,\varphi(z)|\gg 1$. Therefore, in this case we can use the estimates (2.9), (2.10) and (2.18) to obtain 
  $$|{\rm Ai}(\nu^{2/3}\zeta)|^{-1}
 \le C\nu^{1/6}|\zeta|^{1/4}e^{-\frac{2\nu}{3}|{\rm Re}\,\varphi(z)|},
 \eqno{(2.23)}$$
 $$\left|\Phi(\zeta)\right|\le 
 Ce^{-\nu|{\rm Re}\,\varphi(z)|}.
 \eqno{(2.24)}$$
 
 b) $z\in\Theta_2(\delta,\delta_1)$. Then $|\zeta|\ge C_\delta>0$, and by Lemma 2.3 we have $|\arg\zeta(z)|={\cal O}(\delta)$. Hence in this case we can use the expansions (2.5) and 
 (2.6) to obtain 
 $$|{\rm Ai}(\nu^{2/3}\zeta)|^{-1}
 \le C\nu^{1/6}|\zeta|^{1/4}e^{-\frac{2\nu}{3}|{\rm Re}\,\varphi(z)|},
 \eqno{(2.25)}$$
 $$\left|\Phi(\zeta)\right|\le C\nu^{-2}.
 \eqno{(2.26)}$$
 
 c) $z\in\Theta_0(\delta,\delta_1)$. Then we have
 $$\nu^{-1}\le |{\rm Im}\,z|\le|z-1|\le|\zeta|\le 2|z-1|\le 2\delta^2$$
 and by (2.19), $|{\rm Im}\,\zeta|\ge |{\rm Im}\,z|$. Note also that in view of the expansions (2.5) and (2.6), the bounds (2.7) and (2.9)
 hold for all 
$\sigma\in \C\setminus (-\infty,0)$. Using this fact together with (2.19) we obtain in this case
$$|{\rm Ai}(\nu^{2/3}\zeta)|^{-1}
 \le C\nu^{1/3}e^{-\frac{2\nu}{3}|{\rm Re}\,\varphi(z)|},
 \eqno{(2.27)}$$
$$\left|\nu^{-1/3}\zeta^{-1/2}F(\nu^{2/3}\zeta)\right|
\le C+C|\zeta|^{-1/2}(\nu|{\rm Im}\,\zeta|)^{-1}$$ 
$$\le C+Cw(\zeta)^{-1/2}(\nu|{\rm Im}\,z|)^{-1}
 \eqno{(2.28)}$$
and
$$\left|\Phi(\zeta)\right|
\le C+C|\zeta|^{-1/2}(\nu|{\rm Im}\,\zeta|)^{-1}+(4\nu|\zeta|^{3/2})^{-1}$$ 
$$\le C+Cw(\zeta)^{-1/2}(\nu{|\rm Im}\,z|)^{-1}+C\nu^{-1}w(\zeta)^{-3/2}.
 \eqno{(2.29)}$$
It follows from the above bounds that in all three cases we have, for $|{\rm Im}\,z|\gg \nu^{-1}$,
 $$\nu^{-1/3}|\zeta|^{-1/2}|F(\nu^{2/3}\zeta)|\le  Cw(\zeta)^{-1/2},\eqno{(2.30)}$$
 $$|{\rm Ai}(\nu^{2/3}\zeta)|^{-1}
 \le C\nu^{1/3}\langle\zeta\rangle^{1/4}e^{-\frac{2\nu}{3}|{\rm Re}\,\varphi(z)|}.
 \eqno{(2.31)}$$
 In view of (2.8) we also have
 $$|{\rm Ai}(\nu^{2/3}\zeta)|
 \le C\langle\zeta\rangle^{-1/4}e^{\frac{2\nu}{3}|{\rm Re}\,\varphi(z)|}.
 \eqno{(2.32)}$$
 By (2.22) and (2.31), we get, for $|{\rm Im}\,z|\gg \nu^{-1}$, 
 $$|P_3(\zeta)|+|Q_3(\zeta)|\le C_M\nu^{-2M+1}.\eqno{(2.33)}$$
By (2.30) and (2.33), for $|{\rm Im}\,z|\gg \nu^{-1}$ and $\nu$ large enough, we get
 $$\left(1+\frac{1}{|z|}\right)\left|\psi_\nu(\nu z)-\frac{(1-z^2)^{1/2}}{z}\right|$$
$$\le \frac{2}{w(z)}\left|\frac{(1-z^2)^{1/2}}{z}\right|\left(\left|\Phi(\zeta)\right|
 +\widetilde C\nu^{-1}w(z)^{2-\epsilon}w(\zeta)^{-3/2}+\widetilde C\nu^{-2}\right)$$
$$\le 4w(\zeta)^{1/2}w(z)^{-2}\left(\left|\Phi(\zeta)\right|
 +\widetilde C\nu^{-1}w(z)^{2-\epsilon}w(\zeta)^{-3/2}+\widetilde C\nu^{-2}\right)$$ 
 $$\le 4w(\zeta)^{1/2}w(z)^{-2}\left|\Phi(\zeta)\right|
 +4\widetilde C\nu^{-1}w(z)^{-\epsilon}w(\zeta)^{-1}+4\widetilde C\nu^{-2}w(z)^{-2}.$$
 Taking into account that $w(z)\sim 1$, $w(\zeta)\sim 1$, $|{\rm Re}\,\varphi(z)|\ge C|{\rm Im}\,z|$ in case a), 
 $w(\zeta)\sim 1$, $w(z)\sim|z|$, $|z|\gg\nu^{-1}$ in case b), and $w(z)\sim 1$, $w(\zeta)\sim|\zeta|\le 2\delta^2$, $|\zeta|\ge |{\rm Im}\,z|$ in case 
 c), we deduce from the above estimate combined with (2.24), (2.26) and (2.29),
 $$\left(1+\frac{1}{|z|}\right)\left|\psi_\nu(\nu z)-\frac{(1-z^2)^{1/2}}{z}\right|\le
 \left\{
 \begin{array}{lll}
 C_\delta e^{-C_\delta\nu|{\rm Im}\,z|}+C_\delta\nu^{-1},\quad \mbox{in case a)},\\
 C_\delta(\nu|z|)^{-2}+C_{\epsilon,\delta} \nu^{-1+\epsilon},\quad \mbox{in case b)},\\
 C\delta+C_\delta(\nu|{\rm Im}\,z|)^{-1},\quad \mbox{in case c)},
\end{array}
\right.
\eqno{(2.34)}$$
 where the constant $C>0$ is independent of $\delta$. Now,  we can make the LHS of (2.34) less than $(C+1)\delta$ by
 taking $\nu|{\rm Im}\,z|$ and $\nu$ large enough. This implies (2.1) in view of (2.17) after making the change
 $(C+1)\delta\to\delta$.
 
 Given $0<\kappa<1$, define the functions $\varphi_\kappa(z)$ and $\zeta_\kappa(z)$ by 
 $\varphi_\kappa(z)=\varphi(\kappa z)$ and $\zeta_\kappa(z)=\zeta(\kappa z)$.  
 To bound the function $\eta_\nu^{(0)}(\nu z)$, we write it in the form
 $$\eta_\nu^{(0)}(\nu z)=\frac{{\rm Ai}\left(\nu^{2/3}\zeta_\kappa\right)}{{\rm Ai}\left(\nu^{2/3}\zeta\right)}
 \frac{\left(1+Q_1(\zeta_\kappa)+\nu^{-1/3}\zeta_\kappa^{-1/2}F(\nu^{2/3}\zeta_\kappa)Q_2(\zeta_\kappa)+Q_3(\zeta_\kappa)\right)}
 {\left(1+Q_1(\zeta)+\nu^{-1/3}\zeta^{-1/2}F(\nu^{2/3}\zeta)Q_2(\zeta)+Q_3(\zeta)\right)}.$$
 As above, using (2.30)-(2.33), we have, for $\nu\gg 1$, $|{\rm Im}\,z|\gg\nu^{-1}$, 
 $$|\eta_\nu^{(0)}(\nu z)|\le 2\left|\frac{{\rm Ai}\left(\nu^{2/3}\zeta_\kappa\right)}{{\rm Ai}\left(\nu^{2/3}\zeta\right)}\right|
 \le C\nu^{1/3}\left(\frac{\langle\zeta\rangle}{\langle\zeta_\kappa\rangle}\right)^{1/4}e^{-\frac{2\nu}{3}{\rm Re}(\varphi_\kappa(z)-\varphi(z))}$$
  $$\le C\nu^{1/3}\left(\frac{\langle\varphi\rangle}{\langle\varphi_\kappa\rangle}\right)^{1/6}e^{-\frac{2\nu}{3}{\rm Re}(\varphi_\kappa(z)-\varphi(z))}.\eqno{(2.35)}$$
  On the other hand, in view of (2.15), we have the formula
  $$\varphi_\kappa(z)-\varphi(z)=-\int_{\kappa}^1\frac{\varphi_\tau(z)}{d\tau}d\tau=-z\int_{\kappa}^1\varphi'(\tau z)d\tau
  =\int_{\kappa}^1\sqrt{1-(\tau z)^2}\,d\tau.\eqno{(2.36)}$$
  It follows from (2.36) that
  $$\left|\varphi_\kappa(z)-\varphi(z)\right|\le C_1\langle z\rangle$$
  which in turn implies
  $$\frac{\langle\varphi\rangle}{\langle\varphi_\kappa\rangle}\le 1+C_2\frac{\langle z\rangle}{\langle\varphi_\kappa\rangle}\le C_3
  \eqno{(2.37)}$$
  since $\langle\varphi_\kappa\rangle\sim\kappa|z|$ as $|z|\to+\infty$. Set $\Theta_j:=\Theta_j(\delta,\delta_1)$, $j=0,1,2$,
  for some fixed, sufficiently small constants $\delta,\delta_1>0$. It is easy to see that
  $${\rm Re}\,\sqrt{1-(\tau z)^2}\ge 
  \left\{
 \begin{array}{lll}
 C|{\rm Im}\,z|,\quad z\in \Theta_1\cup\Theta_0,\\
 C,\quad z\in \Theta_2,
 \end{array}
 \right.$$
 for all $\kappa\le\tau\le 1$, with a constant $C>0$ independent of $z$ and $\tau$. Hence, by (2.36), 
 $${\rm Re}\,(\varphi_\kappa(z)-\varphi(z))\ge 
  \left\{
 \begin{array}{lll}
 \widetilde C|{\rm Im}\,z|,\quad z\in \Theta_1\cup\Theta_0,\\
 \widetilde C,\quad z\in \Theta_2,
 \end{array}
 \right.
 \eqno{(2.38)}$$
  with a constant $\widetilde C>0$ independent of $z$. By (2.35), (2.37) and (2.38), we conclude
  $$|\eta_\nu^{(0)}(\nu z)|\le
   \left\{
 \begin{array}{lll}
 C'\nu^{1/3}e^{-C\nu|{\rm Im}\,z|},\quad z\in \Theta_1\cup\Theta_0,\\
 C'e^{-C\nu},\quad z\in \Theta_2,
 \end{array}
 \right.
 \eqno{(2.39)}$$
  with constants $C, C'>0$ independent of $z$ and $\nu$. In particular, (2.39) implies
  $$\left(1+\frac{1}{|z|}\right)^2|\eta_\nu^{(0)}(\nu z)|\le C''(\nu|z|)^{1/3}e^{-C\nu|{\rm Im}\,z|}$$
  for all $z$ such that $\nu^{-1}\ll|{\rm Im}\,z|\ll{\rm Re}\,z$, which is the desired bound (2.2) for $|\eta^{(0)}_{\nu}(\lambda)|$.
  In view of the formula $\eta^{(1)}_{\nu}(\lambda)=\psi_\nu(\kappa\lambda)\eta^{(0)}_{\nu}(\lambda)$, the bound for $|\eta_{\nu}^{(1)}(\lambda)|$ follows from that one for 
 $|\eta_{\nu}^{(0)}(\lambda)|$ and the fact that (2.1) implies the bound $|\psi_\nu(\kappa\lambda)|\le C(1+\nu/|\lambda|)$.
 We can bound $|\eta^{(2)}_{\nu}(\lambda)|$ similarly because of the formula
 $$\eta^{(2)}_{\nu}(\lambda)=\eta^{(0)}_{\nu}(\lambda)\left(\left(\frac{\nu}{\kappa\lambda}\right)^2-1-
 (\kappa\lambda)^{-1}\psi_\nu(\kappa\lambda)\right).$$
 
  \eproof

\section{Some properties of the interior Dirichlet-to-Neumann map}

Let $\Omega=\{x\in \R^d:\,|x|\le 1\}$, $\Gamma=\partial\Omega$, and let $\lambda\in\C$ with $1\ll|{\rm Im}\,\lambda|\ll
{\rm Re}\,\lambda$. 
Given a function $f\in H^{s+1}(\Gamma)$, let $u$ solve the equation
$$\left\{
\begin{array}{lll}
 \left(\Delta+\lambda^2\right)u=0&\mbox{in}& \Omega,\\
 u=f&\mbox{on}&\Gamma,
\end{array}
\right.
\eqno{(3.1)}
$$
where $\Delta$ is the negative Euclidean Laplacian. 
We define the interior Dirichlet-to-Neumann (DN) map 
$${\cal N}_0(\lambda):H^{s+1}(\Gamma)\to H^s(\Gamma)$$
 by
$${\cal N}_0(\lambda)f:=\lambda^{-1}\partial_\nu u|_{\Gamma}$$
  $\nu$ being the unit inner normal to $\Gamma$. 
Let $\Delta_\Gamma$ be the negative Laplace-Beltrami operator on the boundary
  $\Gamma$ equipped with the Riemannian metric induced by the Euclidean one.
  In what follows we will denote by $H^1_{sc}(\Gamma)$ the Sobolev space equipped with the semi-classical norm
  $$\|f\|_{H^1_{sc}(\Gamma)}=\|(I-|\lambda|^{-2}\Delta_\Gamma)^{1/2}f\|_{L^2(\Gamma)}$$
   where $I$ denotes the identity.
  For $\sigma\ge 0$, set $$\rho_0(\sigma)=\sqrt{\left(\sigma+\left(\frac{d-2}{2}\right)^2\right) \lambda^{-2}-1}\quad\mbox{with}\quad
  {\rm Re}\,\rho_0>0.$$
 
 \begin{Theorem} For every $0<\delta\ll 1$, independent of $\lambda$, there are positive constants $C_\delta$, $\widetilde C_\delta$
 and $\delta_1=\delta_1(\delta)$ 
 such that for ${\rm Re}\,\lambda\ge \widetilde C_\delta$, $C_\delta\le|{\rm Im}\,\lambda|\le\delta_1{\rm Re}\,\lambda$, we have the estimate
 $$\left\|{\cal N}_0(\lambda)+\rho_0(-\Delta_\Gamma)-
 \frac{d-2}{2\lambda}I\right\|_{L^2(\Gamma)\to H^1_{sc}(\Gamma)}\le \delta.\eqno{(3.2)}$$
\end{Theorem}
 
 {\it Proof.} We will express the DN map in terms of the Bessel functions. If $r=|x|$ is the radial variable, we have
 $$r^{\frac{d-1}{2}}\Delta r^{-\frac{d-1}{2}}=\partial_r^2+\frac{\Delta_\Gamma-(d-1)(d-3)/4}{r^2}.\eqno{(3.3)}$$
 Let $\{\mu_j^2\}$ be the eigenvalues of $-\Delta_\Gamma$ repeated with their multiplicities and let $\{e_j\}$, $\|e_j\|=1$, be the corresponding eigenfunctions, that
 is, $-\Delta_\Gamma e_j=\mu_j^2 e_j$. Denote by $\langle\cdot , \cdot\rangle$ and $\|\cdot\|$ the scalar product and the norm in $L^2(\Gamma).$ If the functions $u$ and $f$ satisfy equation (3.1), we write
 $$f=\sum_j f_je_j,\quad f_j=\langle f,e_j\rangle,\quad \|f\|^2=\sum_j |f_j|^2,$$
 $$u=\sum_j u_j (r)e_j,\quad u_j(r)=\langle u(r, \cdot),e_j(\cdot)\rangle.$$
 In view of (3.3), $w_j(r)=r^{\frac{d-1}{2}}u_j(r)$ and $f_j$ satisfy the equation
 $$\left\{
\begin{array}{lll}
 \left(\partial_r^2-(\nu_j^2-1/4)r^{-2}+\lambda^2\right)w_j=0&\mbox{in}& (0,1),\\
 w_j=f_j&\mbox{at}&r=1,
\end{array}
\right.
\eqno{(3.4)}
$$
 where
 $$\nu_j=\sqrt{\mu_j^2+\left(\frac{d-2}{2}\right)^2}.$$
 The solution of (3.4) is given by the formula
 $$w_j(r)=\frac{b_{\nu_j}(r\lambda)}{b_{\nu_j}(\lambda)}f_j=r^{1/2}\frac{J_{\nu_j}(r\lambda)}{J_{\nu_j}(\lambda)}f_j$$
 where $b_\nu$ and $J_\nu$ are the functions introduced in the previous section. Hence
 $$u_j(r)=r^{-\frac{d-2}{2}}\frac{J_{\nu_j}(r\lambda)}{J_{\nu_j}(\lambda)}f_j.$$
 Since $\partial_\nu u|_\Gamma=-\partial_r u|_{r=1}$, we have
 $${\cal N}_0(\lambda)f=-\sum_j\lambda^{-1}\partial_ru_j|_{r=1}f_j=\sum_j\left(-\psi_{\nu_j}(\lambda)+\frac{d-2}{2\lambda}\right)f_j$$
 where $\psi_{\nu}(\lambda)=J'_{\nu}(\lambda)/J_{\nu}(\lambda)$. This implies 
 $$\left\|\left(I-|\lambda|^{-2}\Delta_\Gamma\right)^{1/2}\left({\cal N}_0(\lambda)+
 \rho_0(-\Delta_\Gamma)-\frac{d-2}{2\lambda}I\right)f\right\|_{L^2(\Gamma)}^2$$
 $$=
 \sum_j\left(1+|\lambda|^{-2}\mu_j^2\right)\left|\psi_{\nu_j}(\lambda)-\rho_0(\mu_j^2)
 \right|^2|f_j|^2$$
 $$\le \sup_{\nu\ge 0}\left(1+|\lambda|^{-2}\nu^2\right)\left|\psi_{\nu}(\lambda)-\rho(\nu/\lambda)\right|^2\|f\|^2
 \eqno{(3.5)}$$ 
  where the function $\rho$ is as in the previous section. Now (3.2) follows from (3.5) and Theorem 2.1.
  \eproof
  
  Let $0<\kappa_1<\kappa_2<1$ be constants and let $\phi(r)\in C_0^\infty([\kappa_1,\kappa_2])$. Then the function $\chi(x)=\phi(|x|)$
  vanishes near $\Gamma$. Given an integer $k\ge0$, denote by $H^k_{r,sc}(\Omega)$ the space equipped with the semi-classical norm 
  $$\|u\|_{H^k_{r,sc}(\Omega)}=\sum_{\ell=0}^k
  |\lambda|^{-\ell}\|\partial_r^\ell u\|_{L^2(\Omega)}$$
 where $r=|x|$ is the radial variable. 
  It is easy to see that the estimate (2.2) implies the following
  
  \begin{lemma} There exist positive constants $C$ and $\widetilde C$ such that the solution $u$ of equation $(3.1)$
  satisfies the estimate
  $$\|\chi u\|_{H^2_{r,sc}(\Omega)}\le \widetilde C|\lambda|^{1/3}e^{-C|{\rm Im}\,\lambda|}\|f\|_{L^2(\Gamma)}.\eqno{(3.6)}$$ 
  \end{lemma}
  
  We will now study the DN map in a more general situation. Let $c(x),\:n(x) \in C^\infty(\overline\Omega)$ be strictly positive functions
  and define the DN map associated to these functions by
  $${\cal N}(\lambda)f:=\lambda^{-1}\partial_\nu u|_{\Gamma}$$
  where $u$ is the solution to the equation
  $$\left\{
\begin{array}{lll}
 \left(\nabla c(x)\nabla+n(x)\lambda^2\right)u=0&\mbox{in}& \Omega,\\
 u=f&\mbox{on}&\Gamma.
\end{array}
\right.
\eqno{(3.7)}
$$
We suppose that there exist a constant $0<\delta_0\ll 1$ and positive constants $\widetilde c$ and $\widetilde n$ such that
$c(x)=\widetilde c$, $n(x)=\widetilde n$ in $\Omega(\delta_0)$. Set
$$\widetilde \rho(\sigma)=\sqrt{\left(\sigma+\left(\frac{d-2}{2}\right)^2\right) \lambda^{-2}-\widetilde n/\widetilde c}\quad\mbox{with}\quad
  {\rm Re}\,\widetilde\rho>0.$$
  
 \begin{Theorem} For every $0<\delta\ll 1$, independent of $\lambda$, there are positive constants $C_\delta$, $\widetilde C_\delta$
 and $\delta_1=\delta_1(\delta)$ 
 such that for ${\rm Re}\,\lambda\ge \widetilde C_\delta$, $C_\delta\log|\lambda|\le|{\rm Im}\,\lambda|\le\delta_1{\rm Re}\,\lambda$, we have the estimate
 $$\left\|{\cal N}(\lambda)+\widetilde\rho(-\Delta_\Gamma)-
 \frac{d-2}{2\lambda}I\right\|_{L^2(\Gamma)\to H^1_{sc}(\Gamma)}\le \delta.\eqno{(3.8)}$$
 \end{Theorem}
 
 {\it Proof.} We will compare ${\cal N}(\lambda)$ with the DN map $\widetilde {\cal N}(\lambda)$ defined by
 $$\widetilde {\cal N}(\lambda)f:=\lambda^{-1}\partial_\nu u|_{\Gamma}$$
  where $u$ is the solution of the equation
  $$\left\{
\begin{array}{lll}
 \left(\widetilde c\Delta+\widetilde n\lambda^2\right)u=0&\mbox{in}& \Omega,\\
 u=f&\mbox{on}&\Gamma.
\end{array}
\right.
\eqno{(3.9)}
$$
 Clearly, we have
 $$\widetilde {\cal N}(\lambda)=\left(\frac{\widetilde n}{\widetilde c}\right)^{-1/2}{\cal N}_0\left(\lambda\left(\frac{\widetilde n}{\widetilde c}\right)^{1/2}\right).$$ 
 In other words, the estimate (3.2) holds true with ${\cal N}_0$ and $\rho_0$ replaced by $\widetilde {\cal N}$ and $\widetilde \rho$, respectively.
 Therefore, one can easily see that Theorem 3.3 follows from Theorem 3.1 and the following
 
 \begin{lemma} There exist positive constants $C$ and $\widetilde C$ such that we have the estimate
  $$\|{\cal N}(\lambda)-\widetilde {\cal N}(\lambda)\|_{L^2(\Gamma)\to H^1_{sc}(\Gamma)}\le \widetilde C|\lambda|^3e^{-C|{\rm Im}\,\lambda|}.\eqno{(3.10)}$$ 
  \end{lemma}
 
 {\it Proof.} Denote by $G_D$ and $\widetilde G_D$ the Dirichlet self-adjoint realizations of the operators $-n^{-1}\nabla c\nabla$ 
 and $-\widetilde n^{-1}\widetilde c\Delta$ on the Hilbert spaces
 $L^2(\Omega,n(x)dx)$ and $L^2(\Omega,dx)$, respectively. Let $\chi_1$ be a smooth function depending only on the radial variable such that
 $\chi_1=1$ in $\Omega(\delta_0/3)$, $\chi_1=0$ in $\Omega\setminus\Omega(\delta_0/2)$. Let $u_1$ be the solution to (3.7) and 
 $u_2$ the solution to (3.9), $u_1=u_2=f$ on $\Gamma$. We have $u_1-\chi_1u_2=0$ on $\Gamma$ and 
 $$U:=\left(n^{-1}\nabla c\nabla+\lambda^2\right)\chi_1u_2=\left(\widetilde n^{-1}\widetilde c\Delta+\lambda^2\right)\chi_1u_2$$
 $$=-\widetilde n^{-1}\widetilde c(\widetilde G_D-\lambda^2)^{-1}(n^{-1}\widetilde c\Delta+\lambda^2)[\Delta,\chi_1]u_2$$
 $$=-\left(\widetilde n^{-1}\widetilde c\right)^2(\widetilde G_D-\lambda^2)^{-1}V$$
 where $V=[\Delta,[\Delta,\chi_1]]u_2$.
  Hence
 $$u_1-\chi_1u_2=(G_D-\lambda^2)^{-1}U$$
 which implies 
 $${\cal N}(\lambda)f-\widetilde {\cal N}(\lambda)f=\lambda^{-1}\gamma\partial_\nu(G_D-\lambda^2)^{-1}U
 \eqno{(3.11)}$$
 where $\gamma$ denotes the restriction on $\Gamma$. 
 By (3.11) we obtain
 $$\|{\cal N}(\lambda)f-\widetilde {\cal N}(\lambda)f\|_{H^1_{sc}(\Gamma)}$$ $$\le {\cal O}(|\lambda|^{-1})\left\|\gamma\partial_\nu(G_D-\lambda^2)^{-1}
 \right\|_{H^2(\Omega)\to H^1_{sc}(\Gamma)}\left\|(\widetilde G_D-\lambda^2)^{-1}
 \right\|_{L^2(\Omega)\to H^2(\Omega)}\|V\|_{L^2(\Omega)}
 \eqno{(3.12)}$$
 where the Sobolev space $H^{2}(\Omega)$ is equipped with the usual norm. We will use now the fact that the norm in $H^1_{sc}(\Gamma)$
 is bounded from above by the usual norm in $H^1(\Gamma)$. 
Thus, by the trace theorem and the coercivity of $G_D$ and $\widetilde G_D$ we have
 $$ \left\|\gamma\partial_\nu(G_D-\lambda^2)^{-1}
 \right\|_{H^2(\Omega)\to H^1_{sc}(\Gamma)}\le {\cal O}(1)\left\|\gamma\partial_\nu(G_D-\lambda^2)^{-1}
 \right\|_{H^2(\Omega)\to H^1(\Gamma)}$$
 $$\le{\cal O}(1)\left\|(G_D-\lambda^2)^{-1}
 \right\|_{H^2(\Omega)\to H^3(\Omega)}$$
 $$\le{\cal O}(1)\left\|(G_D-\lambda^2)^{-1}
 \right\|_{L^2(\Omega)\to H^1(\Omega)}\le {\cal O}(1)\eqno{(3.13)}$$
 and
 $$\left\|(\widetilde G_D-\lambda^2)^{-1}
 \right\|_{L^2(\Omega)\to H^2(\Omega)}$$ $$\le {\cal O}(1)+{\cal O}(|\lambda|^{2})
 \left\|(\widetilde G_D-\lambda^2)^{-1}
 \right\|_{L^2(\Omega)\to L^2(\Omega)}\le {\cal O}(|\lambda|).\eqno{(3.14)}$$
 On the other hand, it is easy to see that the function $V$ is of the form 
 $$V=\sum_{\ell=0}^2a_\ell(r)\partial_r^\ell (\chi_2u_2)$$
 where $\chi_2$ is a smooth function depending only on the radial variable such that
 $\chi_2=1$ on supp$[\Delta,\chi_1]$, $\chi_2=0$ in $\Omega(\delta_0/4)$. 
 Hence, by Lemma 3.2,
 $$\|V\|_{L^2(\Omega)}\le {\cal O}(|\lambda|^{2})\|\chi_2u_2\|_{H^2_{r,sc}(\Omega)}
 \le {\cal O}(|\lambda|^{7/3})e^{-C|{\rm Im}\,\lambda|}\|f\|_{L^2(\Gamma)}\eqno{(3.15)}$$ 
 with a new constant $C>0$. 
 Now (3.10) follows from (3.12)-(3.15).
 \eproof
 
\section{Eigenvalue-free regions}

In this section we derive Theorem 1.1 from Theorems 3.1 and 3.3. Let $c_j(x),n_j(x)\in C^\infty(\overline\Omega)$, $j=1,2$, be strictly positive functions
such that $c_j(x)=\widetilde c_j$, $n_j(x)=\widetilde n_j$ in $\Omega(\delta_0)$, where $\widetilde c_j$, $\widetilde n_j$ are
positive constants satisfying either the condition
$$\widetilde c_1=\widetilde c_2, \quad \widetilde n_1\neq\widetilde n_2,\eqno{(4.1)}$$
or the condition
$$(\widetilde c_1-\widetilde c_2)(\widetilde c_1\widetilde n_1-\widetilde c_2\widetilde n_2)<0.\eqno{(4.2)}$$
Denote by ${\cal N}_j(\lambda)$ the DN map associated to the pair $(c_j,n_j)$ defined in Section 3 and introduce the operator
$$T(\lambda)=\widetilde c_1{\cal N}_1(\lambda)-\widetilde c_2{\cal N}_2(\lambda).$$
Clearly, to prove Theorem 1.1 one has to show that, under the conditions (4.1) or (4.2), $T(\lambda)f=0$ implies $f=0$
for $\lambda\in\Lambda_\ell$, $\ell=1,2$, where 
$$\Lambda_1=\{\lambda\in \C:\,{\rm Re}\,\lambda\gg 1,\,1\ll|{\rm Im}\,\lambda|\ll 
{\rm Re}\,\lambda\}$$
 when the functions $c_j$, $n_j$ are constants in $\Omega$, 
$$\Lambda_2=\{\lambda\in \C:\,{\rm Re}\,\lambda\gg 1,\,\log({\rm Re}\,\lambda)\ll|{\rm Im}\,\lambda|\ll 
{\rm Re}\,\lambda\}$$
  when the functions $c_j$, $n_j$ are constants in $\Omega(\delta_0)$, only. Denote by $\widetilde\rho_j$, $j=1,2$, the functions obtained by replacing the pair $(c,n)$ by 
$(c_j,n_j)$ in the definition of the function $\widetilde\rho$ introduced in Section 3. 
If $T(\lambda)f=0$, $\lambda\in\Lambda_\ell$, $\ell=1,2$, by Theorems 3.1 and 3.3, respectively, we have for all  $\delta>0$,
$$\left\|(1-|\lambda|^{-2}\Delta_\Gamma)^{1/2}\left(\widetilde\rho_1(-\Delta_\Gamma)-\widetilde\rho_2(-\Delta_\Gamma)\right)f\right\|_{L^2(\Gamma)}
\le \delta\|f\|_{L^2(\Gamma)}\eqno{(4.3)}$$
if (4.1) holds, and
$$\left\|(\widetilde c_1\widetilde\rho_1(-\Delta_\Gamma)-\widetilde c_2\widetilde\rho_2(-\Delta_\Gamma))f\right\|_{L^2(\Gamma)}
\le \delta\|f\|_{L^2(\Gamma)}\eqno{(4.4)}$$
if (4.2) holds. On the other hand, we have
$$g(\sigma):=\widetilde c_1\widetilde\rho_1(\sigma)-\widetilde c_2\widetilde\rho_2(\sigma)=\frac{(\widetilde c_1^2-\widetilde c_2^2)\left(\sigma + \left(\frac{d-2}{2}\right)^2\right)\lambda^{-2}
-(\widetilde c_1\widetilde n_1-\widetilde c_2\widetilde n_2)}{\widetilde c_1\widetilde\rho_1+\widetilde c_2\widetilde\rho_2}.$$
Hence, under the above conditions, $g(\sigma)\neq 0$, $\forall\sigma\ge 0$, and we have the bound
$$|g(\sigma)|^{-1}\le C\left\langle\frac{\sigma}{|\lambda|^2}\right\rangle^{k/2}\eqno{(4.5)}$$
where $k=1$ if (4.1) holds and $k=-1$ if (4.2) holds. This implies that the operator 
$$(1-|\lambda|^{-2}\Delta_\Gamma)^{-k/2}g(-\Delta_\Gamma)^{-1}$$ is bounded on $L^2(\Gamma)$ uniformly in $\lambda$.
Therefore, in both cases by (4.3) and (4.4) we conclude
$$\|f\|_{L^2(\Gamma)}\le C\delta\|f\|_{L^2(\Gamma)},\quad\forall\delta>0,\quad\lambda\in\Lambda_\ell,\eqno{(4.6)}$$
with a constant $C>0$ independent of $\delta$. Hence, taking $\delta$ small enough we deduce from (4.6) that $\|f\|=0$, which is
the desired result.
\eproof

\end{document}